\documentclass[draft]{amsart}

\usepackage{amssymb}

\newtheorem{thm}{Theorem}%[section]
\newtheorem{cor}[thm]{Corollary}
\newtheorem{prop}[thm]{Proposition}
\newtheorem{lem}[thm]{Lemma}

\theoremstyle{remark}

\theoremstyle{definition}

\numberwithin{equation}{section}
\numberwithin{thm}{section}

\begin{document}

%%%%%%%%%%%%%%%%%%%%%%%%%
% Subject classification 
%%%%%%%%%%%%%%%%%%%%%%%%%

% Provide an AMS subject classification with one or two primary classification 
% numbers and, optionally, one or more secondary classification numbers. 
% Use the following format:  "Primary 42B25. Secondary 42B60, 20E26"

%\subjclass{Primary 35L05,35Q40. Secondary 35B65}

%%%%%%%%%%NOTE:In Sept.06, 2006%%%%%%%%%%%%%%%%%%%%%%%%
%%%%%%%%%%\subjclass[2000] command is not effective
%%%%%%%%%%or, the version of amsart.cls is old
%%%%%%%%%%So I have written as follows instead:
\renewcommand{\thefootnote}{\fnsymbol{footnote}}
\footnote[0]{2000\textit{ Mathematics Subject Classification}.
 Primary 35L05, 35Q40. Secondary 35B65, 35Q55.}
%%%%%%%%%%%%%%%%%%%%%%%%%%%%%%%%%%%%%%%%%%%%%%%%%%%%%%%%
%35L05 Wave equation
%35B65 Smoothness and regularity of solutions of PDE
%35Q40 Equations from quantum mechanics
%35Q55 NLS-like (nonlinear Schrodinger) equations
%35J10 Schrodinger operator

%%%%%%%%%
% Title
%%%%%%%%%

% Title, in lower case, with no explicit linebreaks (\\).  If the title
% is too long to be used as a running head, add a short version of the
% title in brackets, as in \title[shorttitle]{fulltitle}.

\title[Weighted HLS inequality and Strichartz estimate]{Weighted 
%Hardy-Littlewood-Sobolev
HLS inequalities for radial functions 
and Strichartz estimates for wave and Schr\"odinger equations}

%%%%%%%%%%%%%%%%%%%%%%%%%%%%%%
% Author names and addresses 
%%%%%%%%%%%%%%%%%%%%%%%%%%%%%%

% Provide one separate \author{...} \address{...} \email{....} entry for each
% author, i.e., do not combine multiple authors.  Separate address lines by double
% slashes.  Do not attach footnotes to author  names. (For acknowledgements use
% the "\thanks" construct below.)
%

\author{Kunio Hidano}
\address{Department of Mathematics, 
Faculty of Education, Mie University, 
1577 Kurima-machiya-cho, Tsu, Mie 514-8507, Japan}

\email{hidano@edu.mie-u.ac.jp}

\author{Yuki Kurokawa} 
\address{General Education, Yonago National College of Technology, 
4448 Hikona-cho, Yonago, Tottori 683-8502, Japan} 

\email{kurokawa@yonago-k.ac.jp}

\begin{abstract}
This paper is concerned with derivation 
of the global or local in time Strichartz estimates 
for radially symmetric solutions 
of the free wave equation 
from some Morawetz-type estimates via 
weighted Hardy-Littlewood-Sobolev (HLS) inequalities.  
In the same way we also derive the weighted end-point 
Strichartz estimates with gain of derivatives 
for radially symmetric solutions 
of the free Schr\"odinger equation. 

The proof of the weighted HLS inequality 
for radially symmetric functions 
involves an application of 
the weighted inequality due to Stein and Weiss 
and the Hardy-Littlewood maximal inequality 
in the weighted Lebesgue space due to Muckenhoupt. 
Under radial symmetry 
we get significant gains 
over the usual HLS inequality and Strichartz estimate.
\end{abstract}

\maketitle

%%%%%%%%%%%%%%%%%%%%%%%%%%%%%%%%%%%%%%%%%%%%%%%%%%%%%%%%%%%%%%%%%%%%%%%%%
% end Topmatter
%%%%%%%%%%%%%%%%%%%%%%%%%%%%%%%%%%%%%%%%%%%%%%%%%%%%%%%%%%%%%%%%%%%%%%%%%

%%%%%%%%%%%%%%%%%%%%%%%%%%%%%%%%%%%%%%%%%%%%%%%%%%%%%%%%%%%%%%%%%%%%%%%%%
% body of paper
%%%%%%%%%%%%%%%%%%%%%%%%%%%%%%%%%%%%%%%%%%%%%%%%%%%%%%%%%%%%%%%%%%%%%%%%%
\section{Introduction}
%%%%%%%%%%%%%%%%%%%%%%%%%%%%%%%%%%%%%%%%%%%%%%%%%%%%%%%%%%%
%%%%%%%%%%%%%%%%%%%%%%%%%%%%%%%%%%%%%%%%%%%%%%%%%%%%%%%%%%%%
%In this paper we discuss the implications of the weighted 
%Hardy-Littlewood-Sobolev (HLS, for short) inequalities for 
%radially symmetric functions on 
%the Strichartz estimates for radially symmetric solutions 
%of the free wave equation and the free Schr\"odinger equation. 
%%%%%%%%%%%%%%%%%%%%%%%%%%%%%%%%%%%%%%%%%%%%%%%%%%%%%%%
%In this paper we show the Strichartz estimates 
%for radially symmetric solutions 
%of the free wave equation and the free Schr\"odinger equation 
%by using some space-time $L^2$-estimates 
%for the free solutions and 
%the weighted Hardy-Littlewood-Sobolev (HLS) inequalities 
%for radially symmetric functions.
%%%%%%%%%%%%%%%%%%%%%%%%%%%%%%%%%%%%%%%%%%%%%%%%%%%%%%%%
In this paper we discuss the roles of the weighted 
Hardy-Littlewood-Sobolev (HLS, for short) inequalities for 
radially symmetric functions in the derivation of 
the Strichartz estimates for
% radially symmetric solutions of 
the free wave equation and the free Schr\"odinger equation. 

In the first half of this paper we prove the weighted HLS inequality 
for radially symmetric functions 
on ${\mathbb R}^n$ $(n\geq 2)$ (see (2.2) below). 
The proof proceeds by writing out 
the Riesz potentials in polar coordinates, 
integrating out the angular coordinates, 
and reducing the argument to the one-dimensional setting. 
We then make use of the weighted inequality due to 
Stein and Weiss \cite{StW} and the Hardy-Littlewood maximal inequality 
in the weighted Lebesgue space due to Muckenhoupt \cite{Mu}. 
Naturally, the weighted HLS inequality 
thereby obtained has some similarity with 
the one-dimensional part of 
the weighted inequalities due to Stein and Weiss, 
except that the norm on the right-hand side of (2.2)
involves such a homogeneous weight function as 
$|x|^{-(n-1)((1/p)-(1/q))}$. 
At the cost of the presence of such a singular weight 
function on the right-hand side, 
the weighted radial HLS inequality (2.2) holds 
even for the Riesz potential whose 
kernel has a rather singular form 
$|x|^{-n+\mu}$ with $\mu=\alpha+\beta+(1/p)-(1/q)$. 
(Compare it with the kernel $|x|^{-n+{\tilde\mu}}$, 
${\tilde\mu}=\alpha+\beta+(n/p)-(n/q)$, 
of the Riesz potential in the usual weighted HLS inequality  
(2.15) below.)

In the second half of this paper 
we discuss how the weighted radial HLS inequality is used 
to prove the global (in space and time) or 
local (in time) Strichartz estimate for radially symmetric solutions. 
It is well-known that the range of admissible exponents 
in the global Strichartz estimate for the free wave equation 
can be significantly improved in the radial setting. 
(See Theorem 6.6.2 of Sogge \cite{So}, 
Proposition 4 of Klainerman and Machedon \cite{KM1}, 
Theorem 1.3 of Sterbenz \cite{Ster}, 
and Theorem 4 of Fang and Wang \cite{FW}.) 
Adapting an argument of Vilela \cite{V}, 
we explain how to derive the global radial Strichartz 
estimate in space dimension $n\geq 3$ 
from the generalized Morawetz estimate (see (3.7) below) 
via the weighted radial HLS inequality. 
Our analysis therefore yields another proof of 
Theorem 1.3 of Sterbenz \cite{Ster}. 
As for the local-in-time radial Strichartz estimate 
of the free wave equation 
we extend the space-time 
$L^q$ estimate due to Sogge in space dimension $n=3$ 
into the space-time mixed-norm estimate in space dimension 
$n\geq 2$ (see (5.3) below). 
For that purpose we exploit the local-in-time 
space-time $L^2$-estimate (5.7) below 
by combining it with the weighted radial HLS inequality. 
Such a method does not end with applications 
to the free wave equation. 
Combined with the global (in space and time) estimate of 
the local smoothing property (6.4) below, 
the weighted radial HLS inequality 
is useful in proving the weighted end-point Strichartz estimate 
for radially symmetric solutions to the free 
Sch\"odinger equation (see (6.3) below). 
In the radial setting we observe 
a significant gain of regularity over 
the end-point estimate due to Keel and Tao \cite{KT}. 

%%%%%%%%%%%%From the referee%%%%%%%%%%%%%%%%%%%%
The authors have received a couple of 
very instructive suggestions from the referee. 
One is concerned with the flexibility in our approach. 
The approach to proving the Strichartz estimates used here 
does not rely upon explicit representations or 
parametrices for the solution. 
Therefore it can provide Strichartz-type estimates 
for the large family of equations with defocusing radial potentials. 
Another is concerned with the weighted versions 
of the inhomogeneous Strichartz estimates. 
As was first observed by Kato \cite{Ka} and has been explored by 
Oberlin \cite{Ob}, Harmse \cite{Ha}, Foschi \cite{Fo} and Vilela \cite{V2}, 
the (unweighted) inhomogeneous Strichartz estimates are known 
to hold for the larger range of exponent pairs. 
We now enjoy the approach based upon the celebrated lemma of 
Christ and Kiselev \cite{CK}, 
and the referee has kindly suggested that 
radial weighted analogs thereby obtained 
may turn out to be very useful for certain nonlinear problems. 
Indeed, by virtue of the Christ-Kiselev lemma 
one of the present authors has obtained some radial weighted analogs 
in order to study global existence of small solutions to 
nonlinear wave equations \cite{H4} 
and nonlinear Schr\"odinger equations \cite{H3} 
with radially symmetric data of scale-critical regularity.
%%%%%%%%%%%%%%%%%%%%%%%%%%%%%%%%%%%%%%%%%%%%%%%%

We conclude this section by explaining the notation. 
By $L^p({\mathbb R}^n,\omega(x)dx)$
we mean the Lebesgue space of 
all $\mu$-measurable functions 
$(d\mu(x)=\omega(x)dx)$ on ${\mathbb R}^n$. 
We simply denote 
$L^p({\mathbb R}^n,dx)$ by $L^p({\mathbb R}^n)$. 
The mixed norm $\|u\|_{L^q({\mathbb R};L^p({\mathbb R}^n))}$ 
for functions $u$ on ${\mathbb R}\times{\mathbb R}^n$ 
is defined as 
$$
\|u\|_{L^q({\mathbb R};L^p({\mathbb R}^n))}
=
\biggl(
       \int_{{\mathbb R}}
        \biggl(
               \int_{{\mathbb R}^n}|u(t,x)|^pdx
        \biggr)^{q/p}
       dt
\biggr)^{1/q}
$$
with an obvious modification for $q=\infty$ or $p=\infty$. 
By $p'$ we denote the exponent conjugate to $p$, that is 
$(1/p)+(1/p')=1$. 
The operator $|D_x|^s$ $(s\in{\mathbb R})$ is defined by using the Fourier 
transform ${\mathcal F}$ and the inverse Fourier transform 
${\mathcal F}^{-1}$, as usual. 
We denote 
by ${\dot H}_2^s({\mathbb R}^n)$ the homogeneous Sobolev space 
$|D_x|^{-s}L^2({\mathbb R}^n)$. 
The free evolution operators for the wave equation and 
the Schr\"odinger equation are defined as 
\begin{eqnarray}
& &
(W\varphi)(t,x)=W(t)\varphi(x)=
{\mathcal F}^{-1}e^{it|\xi|}{\mathcal F}\varphi,\\
& &
(S\varphi)(t,x)=S(t)\varphi(x)=
{\mathcal F}^{-1}e^{it|\xi|^2}{\mathcal F}\varphi,
\end{eqnarray}
respectively. 

This paper is organized as follows. 
In the next section we prove the weighted HLS inequality 
for radial functions. 
Section 3 is devoted to the proof of the global-in-time 
Strichartz estimate for radial solutions to the free wave equation. 
In Section 4 we draw our attention to the limiting case of the estimates 
obtained in Section 3.
An adaptation of observations due to 
Agemi \cite{A}, Rammaha \cite{Ra} and Takamura \cite{Ta1} 
shows the failure of such critical estimates.
In Section 5 we are concerned with the local-in-time 
Strichartz estimate for radial solutions to the free wave equation. 
In the final section we revisit the problem of deriving 
the end-point Strichartz estimate for radial solutions 
to the free Schr\"odinger equation from 
the global (in space) estimate of local smoothing property. 
Using the weighted radial HLS inequality, 
we show the weighted end-point Strichartz estimate 
with gain of derivatives for radial free solutions. 
%%%%%%%%%%%%%%%%%%%%%%%%%%%%%%%%%%%%%%%%%%%%%%%%%
%%%%%%%%%%%%%%%%%%%%%%%%%%%%%%%%%%%%%%%%%%%%%%%%
%%%%%%%%%%%%%%%%%%%%%%%%%%%%%%%%%%%%%%%%%%%%%%%%
%%%%%%%%%%%%%%%%%%%%%%%%%%%%%%%%%%%%%%%%%%%%%%%%
\section{Weighted HLS inequality}
We let 
\begin{equation}
(T_\gamma v)(x)
=
\int_{{\mathbb R}^n}\frac{v(y)}{|x-y|^\gamma}dy,
\quad
0<\gamma<n.
\end{equation}
The purpose of this section is to prove the weighted 
Hardy-Littlewood-Sobolev (HLS) inequalities for 
radially symmetric functions. 
We show the following:
\begin{thm}
Suppose $n\geq 2$. Let $p$, $q$, $\alpha$ and $\beta$ 
satisfy $1<p<q<\infty$, $\alpha<1/p'$, $\beta<1/q$ and 
$\alpha+\beta\geq 0$. Set 
$\mu=\alpha+\beta+(1/p)-(1/q)$. 
There exists a constant $C$ 
depending only on 
$n$, $p$, $q$, $\alpha$ and $\beta$, 
and the inequality 
\begin{equation}
\|
|x|^{-\beta}T_{n-\mu}v
\|
_{L^q({\mathbb R}^n)}
\leq
C\|
|x|^{\alpha-(n-1)(1/p-1/q)}v
\|
_{L^p({\mathbb R}^n)}
\end{equation}
holds for radially symmetric 
$v\in
L^p({\mathbb R}^n,
|x|^{p(
       \alpha-(n-1)(1/p-1/q)
      )}dx)$.
\end{thm}

\noindent{\it Remark}. Obviously, 
the number $\mu$ in Theorem 2.1 is strictly positive. 
Moreover, we should note that 
$\mu$ is strictly smaller than one. 
Indeed, by the assumption 
$\alpha<1/p'$, $\beta<1/q$ 
we see $\mu<1/p'+1/q+1/p-1/q=1$.

\vspace{0.3cm}

{\it Proof of Theorem 2.1.} We start with the well-known formula:
\begin{equation}
(T_{n-s}v)(x)
=
\frac{\omega_{n-1}}{r}
\int_0^\infty
\lambda^{n-2}w(\lambda)d\lambda
\int_{|r-\lambda|}^{r+\lambda}
\rho^{-n+s+1}h(\rho,\lambda;r)^{(n-3)/2}d\rho
\end{equation}
$(0<s<n)$ for radially symmetric function 
$v(x)=w(r)$. 
Here and in what follows we use the notation 
$r=|x|=\sqrt{x_1^2+\cdots+x_n^2}$, 
\begin{equation}
h(\rho,\lambda;r)
=
1-
\biggl(
\frac{r^2+\lambda^2-\rho^2}{2\lambda r}
\biggr)^2,
\end{equation}
$\omega_1=2$, and $\omega_n$ $(n=2,3,\dots)$ is the area of 
$S^{n-1}
=
\{\,x\in{\mathbb R}^n\,|\,|x|=1\,\}$. 
For the proof of (2.3) consult, e.g., John \cite{J} on page 8. 
Since the function $h(\rho,\lambda;r)^{(n-3)/2}$ 
causes another singularity in the case of $n=2$, 
let us first study the case $n\geq 3$. 
By virtue of the following proposition 
our argument will be reduced to 
the special case of the weighted estimate of Stein and Weiss. 
(See Lemma 2.3 below.)
\begin{prop}
Suppose $n\geq 3$, $1<q<\infty$ and $0<s<1$. 
The inequality 
\begin{equation}
r^{(n-1)/q}(T_{n-s}v)(x)
\leq
C\int_0^\infty
\frac{\lambda^{(n-1)/q}w(\lambda)}
{|r-\lambda|^{1-s}}d\lambda
\quad
(x\in{\mathbb R}^n)
\end{equation}
holds for radially symmetric, non-negative 
function $v(x)=w(r)$.
\end{prop}
{\it Proof of Proposition 2.2.} 
Set $I_i=I_i(r)$ $(i=1,2)$ as 
\begin{eqnarray}
& &
\frac{1}{r}
\int_0^{r/2}
\lambda^{n-2}w(\lambda)d\lambda
\int_{r-\lambda}^{r+\lambda}
\rho^{-n+s+1}
h(\rho,\lambda;r)^{(n-3)/2}d\rho\\
& &
+
\frac{1}{r}
\int_{r/2}^\infty
\lambda^{n-2}w(\lambda)d\lambda
\int_{|r-\lambda|}^{r+\lambda}
\rho^{-n+s+1}
h(\rho,\lambda;r)^{(n-3)/2}d\rho
%\nonumber\\
%& &
=:I_1(r)+I_2(r).\nonumber
\end{eqnarray}
For the estimate of $I_1$ we note that 
$-1\leq (r^2+\lambda^2-\rho^2)/(2\lambda r)\leq 1$ 
for $|r-\lambda|\leq \rho\leq r+\lambda$, 
which implies $h(\rho,\lambda;r)^{(n-3)/2}\leq 1$ 
by virtue of the assumption $n\geq 3$. 
We therefore obtain
\begin{eqnarray}
& &
I_1(r)
\leq
\frac1r
\int_0^{r/2}
\lambda^{n-2}w(\lambda)d\lambda
\int_{r-\lambda}^{r+\lambda}\rho^{-n+s+1}d\rho\\
& &
\hspace{0.8cm}
\leq
\frac{C}{r^{n-s}}
\int_0^{r/2}\lambda^{n-1}w(\lambda)d\lambda.\nonumber
\end{eqnarray}
The last inequality is due to the fact that 
for $0<s<1$ and $0\leq\lambda\leq r/2$ 
\begin{equation}
\int_{r-\lambda}^{r+\lambda}\rho^{-n+s+1}d\rho
\leq
2\lambda(r-\lambda)^{-n+s+1}
\leq
C\lambda r^{-n+s+1}.
\end{equation}
For the estimate of $I_2$ 
let us first observe 
$h(\rho,\lambda;r)=(\rho^2/\lambda^2)h(\lambda,\rho;r)$. 
Indeed, we see that 
\begin{eqnarray}
& &
h(\rho,\lambda;r)\\
& &
=
\frac{4\lambda^2r^2-(r^2+\lambda^2-\rho^2)^2}{4\lambda^2r^2}
=
\frac{\{\rho^2-(r-\lambda)^2\}\{(r+\lambda)^2-\rho^2\}}{4\lambda^2r^2}
\nonumber\\
& &
=
\frac{(\rho+r-\lambda)(\rho-r+\lambda)(r+\lambda+\rho)(r+\lambda-\rho)}
{4\lambda^2r^2}
\nonumber\\
& &
=
\frac{\{(\rho+r)^2-\lambda^2\}\{\lambda^2-(\rho-r)^2\}}{4\lambda^2r^2}
\nonumber\\
& &
=
\frac{\rho^2}{\lambda^2}
\biggl\{
1-
\biggl(
\frac{r^2+\rho^2-\lambda^2}{2\rho r}
\biggr)^2
\biggr\}
=
\frac{\rho^2}{\lambda^2}h(\lambda,\rho;r),\nonumber
\end{eqnarray}
as desired. Since 
$-1\leq (r^2+\rho^2-\lambda^2)/(2\rho r)\leq 1$ 
for 
$|r-\lambda|\leq \rho\leq r+\lambda$, 
we have 
$h(\rho,\lambda;r)\leq \rho^2/\lambda^2$ and therefore
\begin{eqnarray}
& &
I_2(r)
\leq
\frac1r
\int_{r/2}^\infty\lambda^{n-2}w(\lambda)d\lambda
\int_{|r-\lambda|}^{r+\lambda}
\rho^{-n+s+1}
\biggl(
\frac{\rho^2}{\lambda^2}
\biggr)^{(n-3)/2}d\rho\\
& &
\hspace{0.8cm}
=
\frac1r
\int_{r/2}^\infty\lambda w(\lambda)d\lambda
\int_{|r-\lambda|}^{r+\lambda}
\rho^{-2+s}
d\rho.\nonumber
\end{eqnarray}
Keeping the assumption 
$0<s<1$ in mind, 
we proceed as 
\begin{eqnarray}
& &
\int_{|r-\lambda|}^{r+\lambda}
\rho^{-2+s}
d\rho
=
\frac{1}{(1-s)|r-\lambda|^{1-s}}
\biggl\{
1-
\biggl(
\frac{|r-\lambda|}{r+\lambda}
\biggr)^{1-s}
\biggr\}\\
& &
\leq
\frac{1}{(1-s)|r-\lambda|^{1-s}}
\biggl(
1-
\frac{|r-\lambda|}{r+\lambda}
\biggr)
=
\frac{2\min\{\lambda,r\}}{(1-s)|r-\lambda|^{1-s}(r+\lambda)}.\nonumber
\end{eqnarray}
Combining (2.10) with (2.11), we get
\begin{equation}
I_2(r)
\leq
\frac{C}{r}\int_{r/2}^\infty
\frac{r}{|r-\lambda|^{1-s}(r+\lambda)}
\lambda w(\lambda)d\lambda
\leq
C\int_{r/2}^\infty\frac{w(\lambda)}{|r-\lambda|^{1-s}}d\lambda.
\end{equation}
Therefore, we have obtained by 
(2.3), (2.6), (2.7) and (2.12) 
\begin{equation}
(T_{n-s}v)(x)
\leq
\frac{C}{r^{n-s}}
\int_0^{r/2}\lambda^{n-1}w(\lambda)d\lambda
+
C\int_{r/2}^\infty
\frac{w(\lambda)}{|r-\lambda|^{1-s}}d\lambda.
\end{equation}
We are in a position to complete the proof of (2.5). 
It follows from (2.13) that 
\begin{eqnarray}
& &
r^{(n-1)/q}
(T_{n-s}v)(x)\\
& &
\leq
Cr^{((n-1)/q)-n+s}\int_0^{r/2}\lambda^{n-1}w(\lambda)d\lambda
+
Cr^{(n-1)/q}\int_{r/2}^\infty
\frac{w(\lambda)}{|r-\lambda|^{1-s}}d\lambda\nonumber\\
& &
\leq
C\int_0^{r/2}
\frac{1}{r^{1-s}}
\biggl(
\frac{\lambda}{r}
\biggr)^{(n-1)(1-(1/q))}\lambda^{(n-1)/q}w(\lambda)d\lambda\nonumber\\
& &
\hspace{0.5cm}
+
C\int_{r/2}^\infty
\frac{\lambda^{(n-1)/q}w(\lambda)}{|r-\lambda|^{1-s}}d\lambda
\leq
C\int_0^\infty
\frac{\lambda^{(n-1)/q}w(\lambda)}{|r-\lambda|^{1-s}}d\lambda
\nonumber
\end{eqnarray}
as desired. The proof of Proposition 2.2 has been finished.$\hfill\square$

Once we have obtained the point-wise (in $x$) estimate (2.5), 
the three or higher dimensional part of 
Theorem 2.1 is an immediate consequence of 
the following lemma due to Stein and Weiss \cite{StW}.
\begin{lem}
Assume 
$n\geq 1$, $0<\gamma<n$, 
$1<p<\infty$, $\alpha<n/p'$, $\beta<n/q$, 
$\alpha+\beta\geq 0$, and 
$1/q=(1/p)+((\gamma+\alpha+\beta)/n)-1$. 
If $p\leq q<\infty$, 
then the inequality
\begin{equation}
\||x|^{-\beta}T_\gamma v\|_{L^q({\mathbb R}^n)}
\leq
C\||x|^\alpha v\|_{L^p({\mathbb R}^n)}
\end{equation}
holds for any 
$v\in L^p({\mathbb R}^n,|x|^{p\alpha}dx)$.
\end{lem}

It is easily seen that 
the three or higher dimensional part of 
Theorem 2.1 is a consequence of 
(2.5) with $s=\mu$ and (2.15) with $n=1$. 
The proof of Theorem 2.1 has been finished for $n\geq 3$.

To show Theorem 2.1 in the case of $n=2$ we set 
$J_i=J_i(r)$ $(i=1,2,3)$ for radially symmetric function 
$v(x)=w(r)$:
\begin{eqnarray}
& &
\hspace{0.3cm}\frac{1}{r}\int_0^{r/2}w(\lambda)d\lambda
\int_{r-\lambda}^{r+\lambda}\rho^{-1+s}h(\rho,\lambda;r)^{-1/2}d\rho\\
& &
+\frac{1}{r}\int_{r/2}^{2r}w(\lambda)d\lambda
\int_{|r-\lambda|}^{r+\lambda}\rho^{-1+s}h(\rho,\lambda;r)^{-1/2}d\rho
\nonumber\\
& &
+\frac{1}{r}\int_{2r}^{\infty}w(\lambda)d\lambda
\int_{\lambda-r}^{\lambda+r}\rho^{-1+s}h(\rho,\lambda;r)^{-1/2}d\rho
=:J_1(r)+J_2(r)+J_3(r).\nonumber
\end{eqnarray}
It is possible to show the counterpart of 
Proposition 2.2 for $J_1$ and $J_3$. 
It is $J_2$ that we must handle quite differently from before. 
Let us begin with the proof of the following:
\begin{prop}
Suppose $1<q<\infty$ and $0<s<1$. The inequality 
\begin{equation}
r^{1/q}J_i(r)
\leq
C\int_0^\infty
\frac{\lambda^{1/q}w(\lambda)}{|r-\lambda|^{1-s}}d\lambda
\quad(i=1,3)
\end{equation}
holds for non-negative $w$. 
\end{prop}
{\it Proof of Proposition 2.4.} We use the property of 
the beta function $B(\cdot,\cdot)$:
\begin{equation}
\int_a^b
\frac{2\rho}{\sqrt{\rho^2-a^2}\sqrt{b^2-\rho^2}}d\rho
=
B\biggl(\frac12,\frac12\biggr)=\pi.
\end{equation}
Observing 
\begin{eqnarray}
& &
\frac1r
\int_{|r-\lambda|}^{r+\lambda}\rho^{-1+s}
h(\rho,\lambda;r)^{-1/2}d\rho\\
& &
=%\leq
\lambda\int_{|r-\lambda|}^{r+\lambda}
\rho^{-2+s}
\frac{2\rho}{\sqrt{\rho^2-(r-\lambda)^2}\sqrt{(r+\lambda)^2-\rho^2}}d\rho
\nonumber\\
& &
\leq
\frac{\lambda}{|r-\lambda|^{2-s}}B\bigg(\frac12,\frac12\biggr),
\nonumber
\end{eqnarray}
we are led to
\begin{equation}
r^{1/q}J_1(r)
\leq
C\int_0^{r/2}
\biggl(
\frac{r^{1/q}\lambda^{1-(1/q)}}{r-\lambda}
\biggr)
\frac{\lambda^{1/q}w(\lambda)}{(r-\lambda)^{1-s}}d\lambda
\end{equation}
and 
\begin{equation}
r^{1/q}J_3(r)
\leq
C\int_{2r}^{\infty}
\biggl(
\frac{r^{1/q}\lambda^{1-(1/q)}}{\lambda-r}
\biggr)
\frac{\lambda^{1/q}w(\lambda)}{(\lambda-r)^{1-s}}d\lambda.
\end{equation}
Since $r^{1/q}\lambda^{1-(1/q)}/(r-\lambda)\leq C$ 
for $0\leq\lambda\leq r/2$ and 
$r^{1/q}\lambda^{1-(1/q)}/(\lambda-r)\leq C$ 
for $2r\leq \lambda$, 
the inequality (2.17) is a consequence of 
(2.20)--(2.21). 
We have finished the proof of Proposition 2.4.$\hfill\square$

It remains to show the bound for $J_2$. 
\begin{prop}
Let 
$p$, $q$, $\alpha$, $\beta$ and $\mu$ be the same as in 
Theorem $2.1$. The inequality 
\begin{eqnarray}
& &
\biggl\|
r^{(1/q)-1}\int_{r/2}^{2r}w(\lambda)d\lambda
\int_{|r-\lambda|}^{r+\lambda}\rho^{-1+\mu}
h(\rho,\lambda;r)^{-1/2}d\rho
\biggr\|
_{L^q((0,\infty),r^{-q\beta}dr)}\\
& &
\leq
C\|r^{1/q}w\|
_{L^p((0,\infty),r^{p\alpha}dr)}\nonumber
\end{eqnarray}
holds.
\end{prop}
{\it Proof of Proposition 2.5.} Without loss of generality 
we may assume that $w$ is non-negative. 
Identifying the dual space of 
$L^q((0,\infty),r^{-q\beta}dr)$ with 
$L^{q'}((0,\infty),r^{q'\beta}dr)$ and 
reversing the order integration twice, we have
\begin{eqnarray}
& &
\biggl\|
r^{(1/q)-1}\int_{r/2}^{2r}w(\lambda)d\lambda
\int_{|r-\lambda|}^{r+\lambda}
\rho^{-1+\mu}h(\rho,\lambda;r)^{-1/2}d\rho
\biggr\|
_{L^q((0,\infty),r^{-q\beta}dr)}\\
& &
=\sup
\int_0^\infty
r^{(1/q)-1}g(r)dr\int_{r/2}^{2r}w(\lambda)d\lambda
\int_{|r-\lambda|}^{r+\lambda}
\rho^{-1+\mu}h(\rho,\lambda;r)^{-1/2}d\rho\nonumber\\
& &
=
\sup
\int_0^\infty
w(\lambda)d\lambda
\int_{\lambda/2}^{2\lambda}r^{(1/q)-1}g(r)dr
\int_{|r-\lambda|}^{r+\lambda}
\rho^{-1+\mu}h(\rho,\lambda;r)^{-1/2}d\rho\nonumber\\
& &
=
\sup
\biggl(
\int_0^\infty w(\lambda)d\lambda
\int_0^{\lambda/2}\rho^{-1+\mu}d\rho
\int_{\lambda-\rho}^{\lambda+\rho}
r^{(1/q)-1}h(\rho,\lambda;r)^{-1/2}g(r)dr\nonumber\\
& &
\hspace{1cm}
+
\int_0^\infty w(\lambda)d\lambda
\int_{\lambda/2}^{\lambda}\rho^{-1+\mu}d\rho
\int_{\lambda/2}^{\lambda+\rho}
r^{(1/q)-1}h(\rho,\lambda;r)^{-1/2}g(r)dr\nonumber\\
& &
\hspace{1cm}
+
\int_0^\infty w(\lambda)d\lambda
\int_{\lambda}^{3\lambda/2}\rho^{-1+\mu}d\rho
\int_{\lambda/2}^{2\lambda}
r^{(1/q)-1}h(\rho,\lambda;r)^{-1/2}g(r)dr\nonumber\\
& &
\hspace{1cm}
+
\int_0^\infty w(\lambda)d\lambda
\int_{3\lambda/2}^{3\lambda}\rho^{-1+\mu}d\rho
\int_{\rho-\lambda}^{2\lambda}
r^{(1/q)-1}h(\rho,\lambda;r)^{-1/2}g(r)dr\biggr)\nonumber\\
& &
=:\sup(L_1+L_2+L_3+L_4).\nonumber
\end{eqnarray}
Here the supremum is taken over all non-negative 
$g\in L^{q'}((0,\infty),r^{q'\beta}dr)$ with 
$\|g\|_{L^{q'}((0,\infty),r^{q'\beta}dr)}=1$. 

In what follows we shall often use the identity
\begin{equation}
r^{-1}h(\rho,\lambda;r)^{-1/2}
=
\frac{2\lambda}
{\sqrt{(\rho-r+\lambda)(\rho+r-\lambda)(r+\lambda-\rho)(r+\lambda+\rho)}}
\end{equation}
as well as the inequality 
$2\lambda\leq r+\lambda+\rho\leq 6\lambda$ 
for $\lambda/2\leq r\leq 2\lambda$, 
$|r-\lambda|\leq\rho\leq r+\lambda$. 
To begin with, we first estimate $L_1$. 
Observing $r+\lambda-\rho\leq(\lambda+\rho)+\lambda-\rho=2\lambda$, 
$r+\lambda-\rho\geq(\lambda-\rho)+\lambda-\rho\geq\lambda$ 
for $\lambda-\rho\leq r\leq\lambda+\rho$ 
and $0\leq\rho\leq\lambda/2$, 
we obtain for $0\leq\rho\leq\lambda/2$
%%%%%%%%%%%%%%%%%%%%%%%%%%%%%%%%%%%%%%%%
%%%%%%%The following is a model%%%%%%%%
%%%%%%%%%%%%%%%%%%%%%%%%%%%%%%%%%%%%%%
\begin{eqnarray}
& &
\int_{\lambda-\rho}^{\lambda+\rho}
r^{(1/q)-1}h(\rho,\lambda;r)^{-1/2}g(r)dr\\
& &
\leq
C\lambda^{1/q}\int_{\lambda-\rho}^{\lambda+\rho}
\frac{1}{\sqrt{(\rho-r+\lambda)(\rho+r-\lambda)}}g(r)dr\nonumber\\
& &
\leq
C\lambda^{1/q}
\biggl(
\int_{\lambda-\rho}^\lambda
\frac{1}{\sqrt{\rho(r-\lambda+\rho)}}g(r)dr
+
\int_{\lambda}^{\lambda+\rho}
\frac{1}{\sqrt{\rho(\lambda+\rho-r)}}g(r)dr\biggr)\nonumber\\
& &
\leq
C\lambda^{1/q}
\biggl(\frac{1}{2\rho}\int_{(\lambda-\rho)-\rho}^{(\lambda-\rho)+\rho}
\biggl|
\frac{\rho}{\eta-(\lambda-\rho)}
\biggr|^{1/2}
g^*(\eta)d\eta\nonumber\\
& &
\hspace{1.7cm}
+
\frac{1}{2\rho}\int_{(\lambda+\rho)-\rho}^{(\lambda+\rho)+\rho}
\biggl|
\frac{\rho}{\lambda+\rho-\eta}
\biggr|^{1/2}
g^*(\eta)d\eta\biggr)\nonumber\\
& &
\leq
C\lambda^{1/q}
\biggl(
\sup_{\sigma>0}
\frac{1}{2\sigma}
\int_{(\lambda-\rho)-\sigma}^{(\lambda-\rho)+\sigma}
\biggl|
\frac{\sigma}{\eta-(\lambda-\rho)}
\biggr|^{1/2}
g^*(\eta)d\eta\nonumber\\
& &
\hspace{1.7cm}
+
\sup_{\sigma>0}
\frac{1}{2\sigma}\int_{(\lambda+\rho)-\sigma}^{(\lambda+\rho)+\sigma}
\biggl|
\frac{\sigma}{\lambda+\rho-\eta}
\biggr|^{1/2}
g^*(\eta)d\eta\biggr)\nonumber\\
& &
=
C\lambda^{1/q}({\mathcal M}_{1/2}g^*)(\lambda-\rho)
+
C\lambda^{1/q}({\mathcal M}_{1/2}g^*)(\lambda+\rho).\nonumber
\end{eqnarray}
Here we have set $g^*(\eta)=g(\eta)$ for $\eta\geq 0$, 
$g^*(\eta)=g(-\eta)$ for $\eta<0$, 
and for $t\in{\mathbb R}$
\begin{equation}
({\mathcal M}_{1/2}f)(t)
=
\sup_{\sigma>0}
\frac{1}{2\sigma}
\int_{t-\sigma}^{t+\sigma}
\biggl|
\frac{\sigma}{\eta-t}
\biggr|^{1/2}
|f(\eta)|d\eta.
\end{equation}
It follows from the observation 
of Lindblad and Sogge that 
the maximal function ${\mathcal M}_{1/2}f$, 
which is a singular variant of the Hardy-Littlewood 
maximal function
\begin{equation}
({\mathcal M}f)(t)
=
\sup_{\sigma>0}
\frac{1}{2\sigma}
\int_{t-\sigma}^{t+\sigma}|f(\eta)|d\eta,
\end{equation}
has the point-wise estimate
\begin{equation}
({\mathcal M}_{1/2}f)(t)
\leq
C({\mathcal M}f)(t)
\end{equation}
(see page 1062 of \cite{LS2}). 
Combining (2.25) with (2.28), 
we see that $L_1$ has the bound such as 
\begin{eqnarray}
& &
L_1\leq
C\int_0^\infty
\lambda^{1/q}w(\lambda)d\lambda
\int_0^{\lambda/2}\rho^{-1+\mu}
({\mathcal M}g^*)(\lambda-\rho)d\rho\\
& &
\hspace{0.7cm}
+
C\int_0^\infty
\lambda^{1/q}w(\lambda)d\lambda
\int_0^{\lambda/2}\rho^{-1+\mu}
({\mathcal M}g^*)(\lambda+\rho)d\rho\nonumber\\
& &
\hspace{0.5cm}
\leq
C\|\lambda^{1/q}w\|
_{L^p((0,\infty),\lambda^{p\alpha}d\lambda)}
\|T_{1-\mu}({\mathcal M}g^*)\|
_{L^{p'}((0,\infty),\lambda^{-p'\alpha}d\lambda)}\nonumber\\
& &
\hspace{0.5cm}
\leq
C\|\lambda^{1/q}w\|_{L^p((0,\infty),\lambda^{p\alpha}d\lambda)}
\|{\mathcal M}g^*\|_{L^{q'}({\mathbb R},|t|^{q'\beta}dt)}.\nonumber
\end{eqnarray}
At the last inequality we have used the one-dimensional part of 
Lemma 2.3. To finish the estimate of $L_1$ we need:
\begin{lem}
Suppose that $1<p<\infty$ and 
$-1<a<p-1$. 
The oparator ${\mathcal M}$ enjoys the boundedness
\begin{equation}
\|{\mathcal M}f\|
_{L^p({\mathbb R},|t|^a dt)}
\leq
C
\|f\|_{L^p({\mathbb R},|t|^a dt)}.
\end{equation}
\end{lem}
The original proof of (2.30) is due to Muckenhoupt \cite{Mu}. 
See also Chapter 5 of Stein \cite{Ste} for further references. 
Before we use Lemma 2.6 
to bound $\|{\mathcal M}g^*\|_{L^{q'}({\mathbb R},|t|^{q'\beta}dt)}$, 
let us see that 
the condition $-1<q'\beta<q'-1$ is satisfied. 
The condition $q'\beta<q'-1$ is equivalent to $\beta<1/q$ 
which is supposed in Theorem 2.1. 
Moreover, 
to see that the condition $-1<q'\beta$ 
is also satisfied, 
we note that the assumption 
$(1/q)-(1/p)-\beta+\mu=\alpha<1/p'$ 
implies 
$\mu<(1/q')+\beta$. 
Since $\mu$ is positive, 
we finally find that 
$-1/q'<\beta$, as desired. 
We may therefore use Lemma 2.6 
to proceed as 
\begin{eqnarray}
& &
L_1\leq
C
\|\lambda^{1/q}w\|
_{L^p((0,\infty),\lambda^{p\alpha}d\lambda)}
\|g^*\|_{L^{q'}({\mathbb R},|t|^{q'\beta}dt)}\\
& &
\hspace{0.5cm}
\leq
C
\|\lambda^{1/q}w\|
_{L^p((0,\infty),\lambda^{p\alpha}d\lambda)}
\|g\|_{L^{q'}((0,\infty),r^{q'\beta}dr)}.\nonumber
\end{eqnarray}
The estimate of $L_1$ has been completed.
%%%%%%%%%%%%%%%%%%%%%%%%%%%%%%%%%%%%%%%%%%%%%

We next consider the estimate of $L_2$. 
Observing 
$r+\lambda-\rho\leq(\rho+\lambda)+\lambda-\rho=2\lambda$, 
$r+\lambda-\rho\geq(\lambda/2)+\lambda-\lambda=\lambda/2$ 
for $\lambda/2\leq r\leq\lambda+\rho$ 
and $\lambda/2\leq\rho\leq\lambda$, 
we obtain for $\lambda/2\leq\rho\leq\lambda$
%%%%%%%%%%%%(2.32)%%%%%%%%%%%%%%%%%%%%%%%%%%%%%%%%
\begin{eqnarray}
& &
\int_{\lambda/2}^{\lambda+\rho}
r^{(1/q)-1}h(\rho,\lambda;r)^{-1/2}g(r)dr\\
& &
\leq
C\lambda^{1/q}
\int_{\lambda/2}^{\lambda+\rho}
\frac{1}{\sqrt{(\rho-r+\lambda)(\rho+r-\lambda)}}g(r)dr\nonumber\\
& &
\leq
C
\lambda^{(1/q)-1}
\biggl(
       \int_{\lambda/2}^\lambda
       \sqrt{
             \frac{\lambda}
             {r-(\lambda-\rho)}
             }g(r)dr
+
       \int_{\lambda}^{\lambda+\rho}
       \sqrt{
             \frac{\lambda}
                  {(\lambda+\rho)-r}
             }g(r)dr
\biggr)\nonumber\\
& &
\leq
C
\lambda^{1/q}
\biggl(
       \frac{1}{2\rho}
       \int_{(\lambda-\rho)-\rho}^{(\lambda-\rho)+\rho}
       \biggl|
              \frac{\rho}{\eta-(\lambda-\rho)}
       \biggr|^{1/2}
       g^*(\eta)d\eta\nonumber\\
& &
\hspace{1.7cm}
+
      \frac{1}
           {2\rho}
      \int_{(\lambda+\rho)-\rho}^{(\lambda+\rho)+\rho}
      \biggl|
             \frac{\rho}
                  {(\lambda+\rho)-\eta}
      \biggr|^{1/2}
      g^*(\eta)d\eta
\biggr)
      \nonumber\\
& &
\leq
C\lambda^{1/q}({\mathcal M}_{1/2}g^*)(\lambda-\rho)
+
C\lambda^{1/q}({\mathcal M}_{1/2}g^*)(\lambda+\rho)\nonumber
\end{eqnarray}
as in (2.25). 
Note that, at the second inequality, 
we have used
$$
\frac{\lambda}{2}
\leq
\rho-\lambda+\lambda
\leq
\rho-r+\lambda
\leq
\lambda-\frac{\lambda}{2}+\lambda=\frac32\lambda
$$
for $\lambda/2\leq r\leq\lambda$ and $\lambda/2\leq\rho\leq\lambda$, 
and 
$$
\frac{\lambda}{2}\leq\rho+\lambda-\lambda
\leq
\rho+r-\lambda
\leq
\lambda+(\lambda+\rho)-\lambda
\leq
2\lambda
$$
for $\lambda\leq r\leq\lambda+\rho$ 
and $\lambda/2\leq\rho\leq\lambda$. 
By virtue of the estimate (2.32) we can obtain 
\begin{equation}
L_2
\leq
C
\|\lambda^{1/q}w\|
_{L^p((0,\infty),\lambda^{p\alpha}d\lambda)}
\|g\|
_{L^{q'}((0,\infty),r^{q'\beta}dr)}
\end{equation}
as in (2.29), (2.31). The estimate of $L_2$ has been completed. 
%%%%%%%%%%%%%%%%%%%%%%%%%%%%%%%%%%%%%%%%%%%%%%%%%%%%%%%%%%%%%%%

Next let us consider the estimate of $L_3$. 
Note that 
$\rho+r-\lambda\leq 5\lambda/2$, 
$\rho+r-\lambda\geq\lambda+r-\lambda\geq \lambda/2$ 
for $\lambda/2\leq r\leq 2\lambda$ 
and $\lambda\leq\rho\leq 3\lambda/2$. 
Using (2.24), we hence have for 
$\lambda\leq\rho\leq 3\lambda/2$
%%%%%%%%%%%%%%%%%%%%%%%%%%%%%%%%%%%%%%%%%%%%
\begin{eqnarray}
& &
\int_{\lambda/2}^{2\lambda}
r^{(1/q)-1}h(\rho,\lambda;r)^{-1/2}g(r)dr\\
& &
\leq
C\lambda^{1/q}
\int_{\lambda/2}^{2\lambda}
\frac{1}{\sqrt{(\rho-r+\lambda)(r+\lambda-\rho)}}g(r)dr\nonumber\\
& &
\leq
C
\lambda^{(1/q)-1}
\biggl(
       \int_{\lambda/2}^\lambda
       \sqrt{
             \frac{\lambda}
             {r-(\rho-\lambda)}
             }g(r)dr
+
       \int_{\lambda}^{2\lambda}
       \sqrt{
             \frac{\lambda}
                  {(\rho+\lambda)-r}
             }g(r)dr
\biggr)\nonumber\\
& &
\leq
C
\lambda^{1/q}
\biggl(
       \frac{1}{2\lambda}
       \int_{(\rho-\lambda)-\lambda}^{(\rho-\lambda)+\lambda}
       \biggl|
              \frac{\lambda}{\eta-(\rho-\lambda)}
       \biggr|^{1/2}
       g^*(\eta)d\eta\nonumber\\
& &
\hspace{1.7cm}
+
      \frac{1}
           {3\lambda}
      \int_{(\rho+\lambda)-(3\lambda/2)}^{(\rho+\lambda)+(3\lambda/2)}
      \biggl|
             \frac{3\lambda/2}
                  {(\rho+\lambda)-\eta}
      \biggr|^{1/2}
      g^*(\eta)d\eta
\biggr)
      \nonumber\\
& &
\leq
C\lambda^{1/q}({\mathcal M}_{1/2}g^*)(\rho-\lambda)
+
C\lambda^{1/q}({\mathcal M}_{1/2}g^*)(\rho+\lambda).\nonumber
\end{eqnarray}
This leads us to the estimate
\begin{equation}
L_3
\leq
C
\|\lambda^{1/q}w\|
_{L^p((0,\infty),\lambda^{p\alpha}d\lambda)}
\|g\|
_{L^{q'}((0,\infty),r^{q'\beta}dr)}
\end{equation}
as before. 
The estimate of $L_3$ has been completed.
%%%%%%%%%%%%%%%%%%%%%%%%%%%%%%%%%%%%%%%%%%%%

It remains to bound $L_4$. 
Note that, for 
$\rho-\lambda\leq r\leq 2\lambda$ 
and 
$3\lambda/2\leq\rho\leq 3\lambda$, 
we have 
$\rho-r+\lambda\leq\rho-(\rho-\lambda)+\lambda=2\lambda$, 
$\rho-r+\lambda\geq(3\lambda/2)-2\lambda+\lambda=\lambda/2$ 
and $\rho+r-\lambda\leq3\lambda+2\lambda-\lambda=4\lambda$, 
$\rho+r-\lambda\geq\rho+(\rho-\lambda)-\lambda\geq\lambda$. 
We therefore obtain for $3\lambda/2\leq\rho\leq 3\lambda$
\begin{eqnarray}
& &
\int_{\rho-\lambda}^{2\lambda}
r^{(1/q)-1}h(\rho,\lambda;r)^{-1/2}g(r)dr\\
& &
\leq
C\lambda^{(1/q)-1}
\int_{\rho-\lambda}^{2\lambda}
\sqrt{\frac{\lambda}{r-(\rho-\lambda)}}g(r)dr\nonumber\\
& &
\leq
C\lambda^{1/q}\frac{1}{4\lambda}
\int_{(\rho-\lambda)-2\lambda}^{(\rho-\lambda)+2\lambda}
\biggl|
\frac{2\lambda}{\eta-(\rho-\lambda)}
\biggr|^{1/2}g^*(\eta)d\eta
\nonumber\\
& &
\leq
C\lambda^{1/q}({\mathcal M}_{1/2}g^*)(\rho-r),\nonumber
\end{eqnarray}
which yields
\begin{equation}
L_4
\leq
C
\|\lambda^{1/q}w\|
_{L^p((0,\infty),\lambda^{p\alpha}d\lambda)}
\|g\|
_{L^{q'}((0,\infty),r^{q'\beta}dr)}
\end{equation}
as before. 
Combining (2.31), (2.33), (2.35), (2.37) 
with (2.23), 
we have shown (2.22). 
The proof of Proposition 2.5 has been finished.
$\hfill\square$

We are in a position to complete the proof of 
Theorem 2.1 for $n=2$. 
This is a direct consequence of 
(2.3), (2.17), (2.15) with $n=1$, 
and (2.22). 
The proof of Theorem 2.1 has been completed for all $n\geq 2$.
$\hfill\square$

{\it Remark}. The inequality (2.2) with $\alpha=\beta=0$ 
is just the one Vilela has used in \cite{V}. 
Vilela has shown the inequality 
by employing some ideas in Stein and Weiss \cite{StW}. 
(See \cite{V} on page 369.) 
Now that we have completed the proof of Theorem 2.1, 
it is obvious that we can show (2.2) for $\alpha=\beta=0$ 
by employing the classical Hardy-Littlewood inequality 
and the Hardy-Littlewood maximal inequality 
in the standard $L^p({\mathbb R}^n)$ space. 
Hence it is also possible 
to show (2.2) without results in \cite{StW}, 
as far as the case $\alpha=\beta=0$ is concerned. 
It is in the case $\alpha\ne 0$ or $\beta\ne 0$ 
that our proof of (2.2) essentially relies upon 
the result of Stein and Weiss \cite{StW}.
%%%%%%%%%%%%%%%%%%%%%%%%%%%%%%%%%%%
%%%%%%%%%%%%%%%%%%%%%%%%%%%%%%%%%%%
%%%%%%Section 3%%%%%%%%%%%%%%
%%%%%%%%%%%%%%%%%%%%%%%%%%%%%%%%%
%%%%%%%%%%%%%%%%%%%%%%%%%%%%%%%
\section{Strichartz estimates for radial solutions}
Adapting an argument of Vilela \cite{V}, 
we explain how the weighted Hardy-Littlewood-Sobolev inequality 
(2.2) is used to prove the Strichartz estimate 
for the free wave equation with radially symmetric data. 
Let us start our consideration with global-in-time estimates. 
Recalling the definition of the operator $W$ (see (1.1)), 
we shall show
\begin{thm}
Suppose $n\geq 3$ and 
$1/2<(n-1)((1/2)-(1/p))<(n-1)/2$. 
There exists a constant $C$ depending on 
$n$ and $p$, 
and the estimate
\begin{eqnarray}
& &
\|W\varphi\|_{L^2({\mathbb R};L^p({\mathbb R}^n))}
\leq
C\||D_x|^s\varphi\|_{L^2({\mathbb R}^n)},\\
& &
\hspace{1.7cm}
\frac12+\frac{n}{p}=\frac{n}{2}-s\nonumber
\end{eqnarray}
holds for radially symmetric $\varphi\in{\dot H}^s_2({\mathbb R}^n)$.
\end{thm}
It should be mentioned that 
Sterbenz has proved (3.1) in a completely different way 
(see Proposition 1.2 of \cite{Ster}). 
As has been done in \cite{Ster}, 
we can actually obtain the following result 
by the interpolation between (3.1) and the energy identity. 
For any integer $n\geq 3$ we define
\begin{eqnarray}
& &
D_n:=
\biggl\{\,(x,y)\in{\mathbb R}^2\,\,\Bigl|\,\,
0<x\leq\frac12,\,0<y\leq\frac12,\\
& &
\hspace{3.1cm}
\frac{n-1}{2}\biggl(\frac12-y\biggr)
<x<(n-1)\biggl(\frac12-y\biggr)\,\biggr\}\nonumber
\end{eqnarray}
and
\begin{equation}
A_n:=
D_n
\cup
\biggl\{\,(x,y)\in{\mathbb R}^2\,\,\Bigl|\,\,
x=0\,\,\mbox{and}\,\,y=\frac12\,\biggr\}.
\end{equation}
\begin{cor}
Suppose $n\geq 3$ and 
$(1/q,1/p)\in A_n$. 
There exists a constant $C$ depending on $n$, $p$, $q$, 
and the estimate
\begin{eqnarray}
& &
\|W\varphi\|_{L^q({\mathbb R};L^p({\mathbb R}^n))}
\leq
C\||D_x|^s\varphi\|_{L^2({\mathbb R}^n)},\\
& &
\hspace{1.7cm}
\frac1q+\frac{n}{p}=\frac{n}{2}-s\nonumber
\end{eqnarray}
holds 
for radially symmetric $\varphi\in{\dot H}^s_2({\mathbb R}^n)$.
\end{cor}
%%%%%%%%%%%%%%%%%%%%%%%%%%%%%%%%%%%%%
%%%%%%%%%%%%%%%%%%%%%%%%%%%%%%%%%%%%%
{\it Remark}. Without the assumption of radial symmetry 
the Strichartz estimate (3.4) holds, 
provided that 
\begin{eqnarray}
& &
n\geq 2,\,0\leq\frac1q\leq\frac12,\,0\leq\frac1p\leq\frac12,
\,\biggl(\frac1q,\frac1p\biggr)\ne(0,0),\,
\frac2q\leq(n-1)\biggl(\frac12-\frac1p\biggr),\\
& &
\hspace{1cm}
\biggl(\frac1q,\frac1p\biggr)\ne\biggl(\frac14,0\biggr)
\,\,\mbox{if}\,\,n=2,\,\,\,\,\,
\biggl(\frac1q,\frac1p\biggr)\ne\biggl(\frac12,0\biggr)
\,\,\mbox{if}\,\,n\geq 3.\nonumber
\end{eqnarray}
See \cite{Str}, \cite{Pe2}, \cite{LS1}, \cite{GV}, \cite{KM2}, 
and \cite{KT} for the proof. 
We note that the condition $2/q\leq(n-1)(1/2-1/p)$ of (3.5) is necessary. 
Otherwise, it is well-known that, using the method of Knapp, 
one can indeed choose a sequence 
$\{\varphi_j\}\subset{\mathcal S}({\mathbb R}^n)$ of 
{\it non-radial} data for which 
the existence of such a uniform constant $C=C(n,p,q)$ 
as in (3.4) is forbidden. 
Keeping in mind that some non-radial solutions yield 
this counterexample, we mention that 
radial symmetry vastly improves on 
the range of the admissible pairs $(1/q,1/p)$. 
Indeed, it has turned out by the works of 
Klainerman and Machedon \cite{KM1}, Sterbenz \cite{Ster}, 
and Fang and Wang \cite{FW} (see also Sogge \cite{So} on page 125) 
that one actually has the Strichartz estimate (3.4) 
under the assumption of radial symmetry 
in the case of 
\begin{equation}
n\geq 2,\,\,\,
0\leq\frac1q\leq\frac12,\,\,\,
0\leq\frac1p\leq\frac12,
\,\,\,
\biggl(\frac1q,\frac1p\biggr)\ne(0,0),\,\,\,
\frac1q<(n-1)\biggl(\frac12-\frac1p\biggr),
\end{equation}
in addition to the obvious case 
$(1/q,1/p)=(0,1/2)$. 

In Section 4 we shall show the condition 
$1/q<(n-1)((1/2)-(1/p))$ of (3.6) is necessary 
for the global-in-time estimate (3.4) to hold 
for radially symmetric data. 
The prime purpose of this section is to explain 
how we can prove Proposition 1.2 of Sterbenz \cite{Ster} 
using the weighted Hardy-Littlewood-Sobolev inequality (2.2). 

\vspace{0.3cm}

{\it Proof of Theorem 3.1.} We use the following result 
which is a generalization of the classical estimate of Morawetz \cite{Mo}. 
\begin{lem}
Suppose $n\geq 2$ and $1/2<\alpha<n/2$. 
There exists a constant $C$ depending on 
$n$ and $\alpha$, 
and the estimate 
\begin{equation}
\||x|^{-\alpha}W\varphi\|
_{L^2({\mathbb R}\times{\mathbb R}^n)}
\leq
C\||D_x|^{\alpha-(1/2)}\varphi\|_{L^2({\mathbb R}^n)}
\end{equation}
holds for $\varphi\in {\dot H}^{\alpha-(1/2)}_2({\mathbb R}^n)$. 
\end{lem}
The proof of (3.7) uses the trace inequality 
in the Fourier space
\begin{equation}
\sup_{\lambda>0}
\lambda^{(n/2)-s}
\int_{S^{n-1}}|{\hat w}(\lambda\omega)|^2d\sigma
\leq
C\||D_\xi|^s{\hat w}\|_{L^2({\mathbb R}^n)}
=
C'\||x|^sw\|_{L^2({\mathbb R}^n)}
\end{equation}
which holds for $1/2<s<n/2$. 
See Ben-Artzi \cite{BA}, 
Ben-Artzi and Klainerman \cite{BA-K}, 
Hoshiro \cite{Ho} for the proof of (3.7) 
via the trace inequality such as (3.8) 
and the duality argument. 
For the proof of (3.8) 
see, e.g., (2.45) of Li and Zhou \cite{LZ} and 
Appendix of Hidano \cite{H}.

We are in a position to complete the proof of 
Theorem 3.1. 
We follow the argument of Vilela \cite{V}. 
Fix any $p$ satisfying 
$1/2<(n-1)((1/2)-(1/p))<(n-1)/2$. 
It follows from Theorem 2.1 with $\alpha=\beta=0$ 
that 
the Sobolev-type inequality 
\begin{equation}
\|v\|_{L^p({\mathbb R}^n)}
\leq
C\||x|^{-(n-1)((1/2)-(1/p))}
|D_x|^{(1/2)-(1/p)}v\|_{L^2({\mathbb R}^n)}
\end{equation}
holds for radially symmetric $v$. 
The estimate (3.1) is an immediate consequence of 
(3.7) and (3.9). 
Indeed, we see that 
\begin{eqnarray}
& &
\|W\varphi\|_{L^2({\mathbb R};L^p({\mathbb R}^n))}\\
& &
\leq
C
\||x|^{-(n-1)((1/2)-(1/p))}|D_x|^{(1/2)-(1/p)}W\varphi\|
_{L^2({\mathbb R}\times{\mathbb R}^n)}\nonumber\\
& &
=
C
\||x|^{-(n-1)((1/2)-(1/p))}W(|D_x|^{(1/2)-(1/p)}\varphi)\|
_{L^2({\mathbb R}\times{\mathbb R}^n)}\nonumber\\
& &
\leq
C\||D_x|^{((n-1)/2)-(n/p)}\varphi\|_{L^2({\mathbb R}^n)}\nonumber
\end{eqnarray}
as desired. 
The proof of Theorem 3.1 has been finished.$\hfill\square$
%%%%%%%%%%%%%%%%%%%%%%%%%%%%%%%%%%%%%%%%%%%%%%%%%%%%%%%%%%%%
%%%%%%%%%%%%%%%%%%%%%%%%%%%%%%%%%%%%%%%%%%%%%%%%%%%%%%%%%%%%
%%%%%%%%%%%Section 4%%%%%%%%%%%%%%%%%%%%%%%%%%%%%%%%%%%%%%%%
%%%%%%%%%%%%%%%%%%%%%%%%%%%%%%%%%%%%%%%%%%%%%%%%%%%%%%%%%%%%
%%%%%%%%%%%%%%%%%%%%%%%%%%%%%%%%%%%%%%%%%%%%%%%%%%%%%%%%%%%%
\section{Failure of the critical estimate}
The problem to be discussed in this section is 
whether the Strichartz estimate (3.4) holds 
under the assumption of radial symmetry of data 
even for the limiting pair 
$(1/q,1/p)\in (0,1/2]\times [0,1/2)$ 
with $1/q=(n-1)((1/2)-(1/p))$. 
If it were true, 
we would enjoy
\begin{equation}
\|W\varphi\|_{L^q({\mathbb R};L^p({\mathbb R}^n))}
\leq
C\||D_x|^{(1/2)-(1/p)}\varphi\|_{L^2({\mathbb R}^n)}
\end{equation}
for radially symmetric data $\varphi$, and 
the estimate (4.1) would imply the estimate
\begin{eqnarray}
& &
\|u\|_{L^q({\mathbb R};L^p({\mathbb R}^n))}\\
& &
\leq
C
\bigl(
\||D_x|^{(1/2)-(1/p)}f\|_{L^2({\mathbb R}^n)}
+
\||D_x|^{-(1/2)-(1/p)}g\|_{L^2({\mathbb R}^n)}
\bigr)
\nonumber
\end{eqnarray}
for the solution $u$ to the wave equation 
$\square u=0$ in ${\mathbb R}\times{\mathbb R}^n$ 
with radially symmetric data $(f,g)$. 
We shall show that the estimate (4.2) is false 
in the limiting case $(1/q,1/p)\in (0,1/2]\times [0,1/2)$ 
with $1/q=(n-1)((1/2)-(1/p))$, though 
${\mathcal S}({\mathbb R}^n)
\subset
{\dot H}^{-(1/2)-(1/p)}_2({\mathbb R}^n)$ ($n\geq 2$). 
%%%%%%%%%%%%%%%%%%%%%%%%%%%%%%%%%%%%%%%%%%%%%
%%%%%%%Kurokawa file%%%%%%%%%%%%%%%%%%%%%%%%%
%%%%%%%%%%%%%%%%%%%%%%%%%%%%%%%%%%%%%%%%%%%%%
The key to such a result is the following 
%%%%%%%%%%%%%%%%%%%%%%%%%%%%%%%%%%%%%%%%%%%%%%%%%%%%%
%%%%%%%%%%%%Kurokawa Lemma%%%%%%%%%%%%%%%%%%%%
%%%%%%%%%%%%%%%%%%%%%%%%%%%%%%%%%%%%%%%%%%%%%%%%%%%%%
\begin{lem}%Lemma 4.1.
Let $n\geq 2$ and $r=|x|$. 
Suppose $g(x)$ is a smooth, non-negative function 
with ${\rm supp}\,g\subset
\{\,x\in{\mathbb R}^n\,|\,|x|\leq R\,\}$ 
for some $R>0$. 
Suppose also that $g$ is a radially symmetric function 
written as $g(x)=\psi(r)$ for an even function 
$\psi\in C_0^\infty({\mathbb R})$. 
Let $u$ be the solution to $\square u=0$ 
in ${\mathbb R}\times{\mathbb R}^n$ 
with data $(0,g)$ at $t=0$. 
There exists a positive constant $\delta$ 
depending only on $n$ such that the estimate 
\begin{equation}
u(t,x)
\geq
\frac{1}{4r^{(n-1)/2}}
\int_{r-t}^{\min\{R,r+t\}}
\lambda^{(n-1)/2}\psi(\lambda)d\lambda
\end{equation}
holds for any $(t,x)$ with 
$R/(1+\delta)\leq r-t\leq R$, $t>0$.
\end{lem}
%%%%%%%%%%%%%%%%%%%%%%%%%%%%%%%%%%%%%%%%%%%%%%%%%%%%%%%
Let us postpone the proof of Lemma 4.1 for the moment 
and see how it can be used to prove
%%%%%%%%%%%%%%%%%%%%%%%%%%%%%%%%%%%%%%%%%%%%%%%%%%%%%%%%
\begin{thm}%Theorem 4.2.
Let $n\geq 2$ and fix the constant $\delta>0$ given by Lemma $4.1$. 
Suppose that $g(x)\geq 0$ is a smooth, radially symmetric 
function with ${\rm supp}\,g\subset\{\,x\in{\mathbb R}^n\,|\,
|x|\leq 1\,\}$ which is written as 
$g(x)=\psi(r)$ for an even function 
$\psi\in C_0^\infty({\mathbb R})$ satisfying the condition that 
the function $\Psi$ defined as
\begin{equation}
\Psi(\rho):=\int_\rho^1
\lambda^{(n-1)/2}\psi(\lambda)d\lambda
\end{equation}
does not vanish identically for $\rho\in(1/(1+\delta),1)$. 

Let $(1/q,1/p)\in (0,1/2]\times [0,1/2)$ satisfy 
$1/q=(n-1)(1/2-1/p)$. 
Then, for the solution to $\square u=0$ with data $(0,g)$ at $t=0$, 
\begin{equation}
\lim_{T\to +\infty}
\|u\|_{L^q((0,T);L^p({\mathbb R}^n))}=+\infty.
\end{equation}
\end{thm}
%%%%%%%%%%%%%%%%%%%%%%%%%%%%%%%%%%%%%%%%%%%%%
{\it Proof of Theorem 4.2.} We separate two cases: 
$(1/q,1/p)\in (0,1/2]\times (0,1/2)$ satisfying 
$1/q=(n-1)(1/2-1/p)$ for $n\geq 2$ 
and $(1/q,1/p)=(1/2,0)$ for $n=2$. 
We start with the former. 
%%%%%%%%%%%%%%%%%%%%%%%%%%%%%%%%%%%%%%%%%%%%%%%%%%%%
Employing (4.3) and writing $u(t,x)=v(t,r)$, 
we have for 
$t\in(\delta/(2(1+\delta)),T)$ 
by the change of variables $\rho=r-t$
\begin{eqnarray}
& &
\int_{t+(1/(1+\delta))}^{t+1}
v^p(t,r)r^{n-1}dr\\
& &
\geq
\frac{1}{4^p}
\int_{t+(1/(1+\delta))}^{t+1}
\biggl(
\frac{1}{r^{(n-1)/2}}
\int_{r-t}^1
\lambda^{(n-1)/2}\psi(\lambda)d\lambda
\biggr)^pr^{n-1}dr\nonumber\\
& &
=
\frac{1}{4^p}
\int_{1/(1+\delta)}^{1}
\frac{1}{(t+\rho)^{(\frac{n-1}{2})p-(n-1)}}
\Psi^p(\rho)
d\rho
\nonumber\\
& &
\geq
\frac{1}{4^p}
\frac{1}{(t+1)^{(\frac{n-1}{2})p-(n-1)}}
\int_{1/(1+\delta)}^{1}
\Psi^p(\rho)
d\rho.\nonumber
\end{eqnarray}
Setting a strictly positive constant $A$ as
$$
A:=
\biggl(
\int_{1/(1+\delta)}^{1}
\Psi^p(\rho)
d\rho
\biggr)^{1/p},
$$
we then find 
\begin{eqnarray}
& &
\|u\|^q_{L^q((0,T);L^p({\mathbb R}^n))}
\geq
\int_{\delta/(2(1+\delta))}^T
\biggl(
\int_{t+(1/(1+\delta))}^{t+1}
v^p(t,r)r^{n-1}dr
\biggr)^{q/p}dt\\
& &
\geq
\frac{A^q}{4^q}\int_{\delta/(2(1+\delta))}^T
\biggl(
\frac{1}{(t+1)^{\frac{n-1}{2}-\frac{n-1}{p}}}
\biggr)^qdt\nonumber\\
& &
=
\frac{A^q}{4^q}\int_{\delta/(2(1+\delta))}^T
\frac{1}{t+1}dt
=
\frac{A^q}{4^q}
\log
\frac{T+1}{\frac{\delta}{2(1+\delta)}+1}\nonumber
\end{eqnarray}
for all $T\geq\delta/(2(1+\delta))$.

It remains to deal with $(1/q,1/p)=(1/2,0)$ for $n=2$. 
We naturally modify the argument in (4.6)--(4.7) as follows. 
Fix a constant $c_0$ satisfying 
$1/(1+\delta)<c_0<1$ so that $\Psi(c_0)>0$. 
We see, noting 
$\|v(t,\cdot)\|_{L^\infty(t+(1/(1+\delta))<r<t+c_0)}
\geq
v(t,t+c_0)$, 
\begin{eqnarray}
& &
\int_{\delta/(2(1+\delta))}^T
\|v(t,\cdot)\|^2_{L^\infty(t+(1/(1+\delta))<r<t+c_0)}dt\\
& &
\geq
\int_{\delta/(2(1+\delta))}^T
\frac{1}{4^2(t+c_0)}
\biggl(
\int_{c_0}^1
\lambda^{1/2}\psi(\lambda)d\lambda
\biggr)^2dt
=
\frac{1}{4^2}
\Psi^2(c_0)
\log
\frac{T+c_0}{\frac{\delta}{2(1+\delta)}+c_0}\nonumber
\end{eqnarray}
for all $T\geq\delta/(2(1+\delta))$. 
We have completed the proof. $\hfill\square$

\vspace{0.5cm}

Using the solution $u$ described in Theorem 4.2, 
we easily obtain the following result by scaling argument. 
%%%%%%%%%%%%%%%%%%%%%%%%%%%%%%%%%%%%%%%%%%%%%%%%%%%%%%%%%
\begin{cor}%Corollary 4.3.
Let $n\geq 2$ and 
$(1/q,1/p)\in (0,1/2]\times [0,1/2)$ satisfy $1/q=(n-1)(1/2-1/p)$. 
Then, for the solution $u_h$ to 
$\square u=0$ with radially symmetric data $(0,h)$ at $t=0$
\begin{eqnarray*}
& &
\sup
\biggl\{
\,
\frac{
      \|u_h\|_{L^q((0,1);L^p({\mathbb R}^n))}
     }
     {\| |D_x|^{-(1/2)-(1/p)}h \|_{L^2({\mathbb R}^n)}
     }\,\,\bigl|\\
& &
\hspace{1.2cm}
h\in{\mathcal S}({\mathbb R}^n)\setminus\{0\}
\mbox{\,and $h$ is radially symmetric}
\,
\biggr\}
=
+\infty.
\end{eqnarray*}

\end{cor}
%%%%%%%%%%%%%%%%%%%%%%%%%%%%%%%%%%%%%%%%%%%%%%%%%%
This shows that 
the estimate (4.2) is false 
even if the global-in-time norm 
is replaced by the local-in-time norm on the left-hand side. 
The proof of Corollary 4.3 is straightforward, 
and therefore we leave it to the reader.

\vspace{0.5cm}

%%%%%%%%%%%%%%%%%%%%%%%%%%%%%%%%%%%%%%%%%%%%%%%%%%%%%%%%
%%%%%%%%%%Kurokawa file%%%%%%%%%%%%%%%%%%%%%%%%%%%%%%%
%%%%%%%%%%%%%%%%%%%%%%%%%%%%%%%%%%%%%%%%%%%%%%%%%%%%%%%
{\it Proof of Lemma 4.1.} We must establish Lemma 4.1. 
The proof is essentially based on Rammaha's way for Lemma 2 of \cite{Ra} 
together with the treatment of fundamental solutions in even 
space dimensions in Agemi \cite{A}, 
which is summarized in Takamura \cite{Ta2}. 
Following \cite{A}, \cite{Ra} and \cite{Ta1}--\cite{Ta2}, 
we show Lemma 4.1. 

By the representations (6a) and (6b) of radial solutions 
in \cite{Ra}, $u$ is expressed as
\begin{equation}
u(t,x)=
\frac{1}{2r^m}
\int_{|r-t|}^{r+t}
\lambda^m\psi(\lambda)
P_{m-1}
\biggl(
\frac{\lambda^2+r^2-t^2}{2r\lambda}
\biggr)
d\lambda,
\end{equation}
if $n=2m+1$, and 
\begin{eqnarray}
& &
u(t,x)\\
& &
=
\frac{2}{\pi r^{m-1}}
\int_0^t
\frac{\rho d\rho}{\sqrt{t^2-\rho^2}}
\int_{|r-\rho|}^{r+\rho}
\frac{\lambda^m\psi(\lambda)}
{\sqrt{G(\lambda,r,\rho)}}
T_{m-1}\left(\frac{\lambda^2+r^2-\rho^2}{2r\lambda}\right)
d\lambda\nonumber\\
& &
=
\frac{2}{\pi r^{m-1}}
\int_{r-t}^{r+t}\lambda^m\psi(\lambda)d\lambda
\int_{|r-\lambda|}^t
\frac{\rho}{\sqrt{G(\rho,r,\lambda)}\sqrt{t^2-\rho^2}}\nonumber\\
& &
\hspace{6cm}
\times
T_{m-1}\left(\frac{\lambda^2+r^2-\rho^2}{2r\lambda}\right)d\rho,\nonumber
\end{eqnarray}
if $n=2m$ and $r>t$, 
%for $(t,x)\in(0,\infty)\times{\mathbb R}^n$,
where 
$$
G(\lambda,r,\rho)=(\lambda^2-(r-\rho)^2)((r+\rho)^2-\lambda^2)
=G(\rho,r,\lambda),
$$
and $P_k$, $T_k$ are the Legendre and Tschebyscheff polynomials, 
respectively, defined by 
\begin{eqnarray*}
& &
P_k(z)=\frac{1}{2^kk!}\frac{d^k}{dz^k}(z^2-1)^k,\\
& &
T_k(z)=\frac{(-1)^k}{(2k-1)!!}(1-z^2)^{1/2}
\frac{d^k}{dz^k}(1-z^2)^{k-1/2}.
\end{eqnarray*}
See also \cite{Ta1} for details. 

As is well-known, $P_k$ and $T_k$ have the properties:
$|P_k(z)|,\,|T_k(z)|\le 1$ $(|z|\le 1)$ and $P_k(1)=T_k(1)=1$ 
for all $k=1,2,\dots$ 
(see Magnus, Oberhettinger and Soni \cite{MOS}, 
p.\,227, p.\,237, pp.\,256--267). 
By these properties together with the continuity of the two functions, 
one can choose a small constant $\delta$ depending on $n$ 
so that 
\begin{equation}
P_{m-1}(z),\quad T_{m-1}(z)\ge\frac{1}{2}
\quad\mbox{for}\quad\frac{1}{1+\delta}\le z\le 1,\quad 
m\in{\mathbb N},
\end{equation}
in the same manner as Takamura did in Lemma 2.5 of \cite{Ta1}.

In what follows we assume that $R/(1+\delta)\le r-t\le R$ with $t>0$. 
Note that the upper limit of the $\lambda$-integrals in 
(4.9) and (4.10) can be replaced with 
$\min\{R,r+t\}$ by virtue of the support property of data. 
We then have 
\begin{eqnarray}
& &
1\ge
\frac{\lambda^2+r^2-\rho^2}{2r\lambda}
\ge
\frac{\lambda^2+r^2-t^2}{2r\lambda}\\
& &
\hspace{0.25cm}
\ge 
\frac{(r-t)^2+r^2-t^2}{2rR}
=
\frac{r-t}{R}
\ge
\frac{1}{1+\delta}\nonumber
\end{eqnarray}
for $r-t\le\lambda\le \min\{R,r+t\}$ and $0\le\rho\le t$. 
It therefore follows from (4.9)--(4.11) that 
\begin{equation}
u(t,x)\geq
\left\{
%%%%%%%%%%%%%%%%%%%%%%%%%%%%%%%%%%%%%%
\begin{array}{l}
\displaystyle{
              \frac{1}{4r^m}
              \int_{r-t}^{\min\{R,r+t\}}
              \lambda^m\psi(\lambda)d\lambda
              \,\,\,\mbox{if}\,\,\,n=2m+1,
             }\\
{}\\
\displaystyle{
              \frac{1}{\pi r^{m-1}}
              \int_{r-t}^{\min\{R,r+t\}}\lambda^m\psi(\lambda)d\lambda
              \int_{|r-\lambda|}^t
              \frac{\rho d\rho}
                   {\sqrt{G(\rho,r,\lambda)}\sqrt{t^2-\rho^2}}
              }\\
{}\\
\hspace{4.5cm}\mbox{if}\,\,\,n=2m,
\end{array}
%%%%%%%%%%%%%%%%%%%%%%%%%%%%%%%%%%%%%%%%%%%
\right.
\end{equation}
provided $\psi\ge 0$ on the support. 
Therefore the odd dimensional case has been proved. 
For the even dimensional case, the $\rho$-integral in (4.13) 
is estimated as follows: 
\begin{eqnarray*}
(\rho\mbox{-integral})
&\ge& 
\frac{1}{2\sqrt{r\lambda}}\int_{|r-\lambda|}^t
\frac{\rho d\rho}{\sqrt{\rho^2-(r-\lambda)^2}\sqrt{t^2-\rho^2}}\\
&=&
\frac{B(2^{-1},2^{-1})}{4\sqrt{r\lambda}}
=\frac{\pi}{4\sqrt{r\lambda}}.
\end{eqnarray*}
Here by $B(\cdot,\cdot)$ we have meant the beta function as in Section 2. 
This completes the proof for the even dimensional case. 

$\hfill\square$

{\it Remark}. During the preparation of this article, 
the authors found that, 
arguing in a way similar to Takamura \cite{Ta1}, 
Jiao and Zhou had already obtained an estimate 
which is a bit less precise than (4.3) 
%essentially the same as (4.3) 
(see Lemma 2 of \cite{JZ}). 
%%%%%%%%%%%%%%%%%%%%%%%%%%%%%%%%%%%%%%%%%
%%%%%%%%%%Section 5%%%%%%%%%%%%%%%%%%%%%%
%%%%%%%%%%%%%%%%%%%%%%%%%%%%%%%%%%%%%%%%%
\section{Local-in-time Strichartz estimates}
For any integer $n\geq 2$ we define
\begin{equation}
\Omega_n:=
\biggl\{
\,(x,y)\in{\mathbb R}^2\,\,\Bigl|\,\,
0<x\leq\frac12,\,
0<y\leq\frac12,\,
x>(n-1)\biggl(\frac12-y\biggr)\,
\biggr\}
\end{equation}
and
\begin{equation}
\Lambda_n:=
\Omega_n\cup
\biggl\{
\,(x,y)\in{\mathbb R}^2\,\,\Bigl|\,\,
x=0\,\,\mbox{and}\,\,y=\frac12\,
\biggr\}
\end{equation}
The main result of this section is the following.
%%%%%%%%%%%%%%%%%%%%%%%%%%%%%%%%%%%%%%%%%%%
%%%%%%%%%%%%Theorem 5.1.%%%%%%%%%%%%%%%%%%%
%%%%%%%%%%%%%%%%%%%%%%%%%%%%%%%%%%%%%%%%%%%
\begin{thm}
Suppose $n\geq 2$ and 
$(1/q,1/p)\in\Lambda_n$. 
Let $T$ be an arbitrary positive number. 
There exists a constant $C$ depending only on 
$n$, $p$ and $q$, and 
the estimate 
\begin{eqnarray}
& &
\|W\varphi\|_{L^q((0,T),L^p({\mathbb R}^n))}
\leq
CT^\theta
\||D_x|^{(1/2)-(1/p)}\varphi\|_{L^2({\mathbb R}^n)},\\
& &
\hspace{1.8cm}
\frac1q+\frac{n}{p}=\theta+\frac{n}{2}
-\biggl(
\frac12-\frac1p
\biggr)\nonumber
\end{eqnarray}
holds for radially symmetric data 
$\varphi\in{\dot H}^{(1/2)-(1/p)}_2({\mathbb R}^n)$.
\end{thm}
Theorem 5.1 is an extension of the intriguing result of Sogge 
(Proposition 6.3 on the page 125 of \cite{So}) 
who proved the estimate (5.3) for 
$n=3$ and 
\begin{eqnarray}
& &
\biggl(
\frac1q,\frac1p
\biggr)
\in
\biggl\{
(x,y)\in{\mathbb R}^2\,\,\Bigl|\,\,
0<x\leq y\leq\frac12,\,
x>2\biggl(\frac12-y\biggr)
\biggr\}\\
& &
\hspace{1.8cm}
\cup
\biggl\{
(x,y)\in{\mathbb R}^2\,\,\Bigl|\,\,
x=0\,\,\mbox{and}\,\,y=\frac12
\biggr\}.\nonumber
\end{eqnarray}
Actually, Sogge himself proved the estimate (5.3) 
for $n=3$, $1/3<1/q=1/p\leq 1/2$. 
By the interpolation 
between his estimate and the energy estimate 
we easily get (5.3) for $n=3$ and 
$(1/q,1/p)$ satisfying (5.4).

We should explain the significance of the local-in-time 
estimate (5.3). 
If the radially symmetric estimate (4.1) were true 
even for the limiting pair 
$(1/q,1/p)\in (0,1/2]\times[0,1/2)$ 
with $1/q=(n-1)((1/2)-(1/p))$, 
our estimate (5.3) would be a trivial consequence 
of (4.1) and the H\"older inequality in time. 
The fact is that 
the estimate (4.1), even if localized in time, is false 
for any limiting pair 
$(1/q,1/p)\in (0,1/2]\times[0,1/2)$ 
with $1/q=(n-1)((1/2)-(1/p))$ 
as we have seen in Section 4, and  
one can get nothing but a coarse estimate 
\begin{eqnarray}
& &
\|W\varphi\|_{L^q((0,T),L^p({\mathbb R}^n))}
\leq
CT^\theta
\||D_x|^{(1/2)-(1/p)+\varepsilon}\varphi\|_{L^2({\mathbb R}^n)},\\
& &
\frac1q+\frac{n}{p}=\theta+\frac{n}{2}
-\biggl(
\frac12-\frac1p+\varepsilon
\biggr),\,\,\varepsilon>0\,\,\mbox{(sufficiently small)}\nonumber
\end{eqnarray}
for any $(1/q,1/p)\in\Omega_n$ with $(1/q,1/p)\in(0,1/2]\times(0,1/2)$ 
by using both the Strichartz estimate (3.4) 
for $(1/q,1/p)$ permitted in (3.6) and 
the H\"older inequality in time. 
As we have just mentioned, 
Sogge proved the sharper estimate (5.3) in the case of 
$n=3$, $1/3<1/q=1/p\leq 1/2$, 
and the key to his proof 
was a clever use of the identity
\begin{eqnarray}
& &
\hspace{0.8cm}
\widehat{d\sigma}(|\xi|)
=
\int_{S^2}e^{-i\omega\cdot\xi}d\sigma
=
4\pi\frac{\sin|\xi|}{|\xi|}\\
& &
\big(\omega\in S^2=\{\,x\in{\mathbb R}^3\,|\,|x|=1\,\},\,\,
d\sigma=d\sigma(\omega)\big).\nonumber
\end{eqnarray}
Though the formula of $\widehat{d\sigma}(|\xi|)$ 
in terms of the Bessel function is 
well-known for $n=2$ or $n\geq 4$, 
the authors do not know whether 
such a formula is useful in proving our estimate (5.3). 
In the rest of this section 
we see how the weighted inequality (2.2) 
is used to prove the local-in-time estimate (5.3).

\vspace{0.3cm}

{\it Proof of Theorem 5.1.} We use the following result. 
\begin{lem}
Suppose $n\geq 1$ and $0\leq\alpha<1/2$. 
Let $T$ be an arbitrary positive number. 
There exists a constant $C$ depending on $n$ and $\alpha$, 
and the estimate 
\begin{equation}
\||x|^{-\alpha}W\varphi\|
_{L^2((0,T)\times{\mathbb R}^n)}
\leq
CT^{(1/2)-\alpha}\|\varphi\|_{L^2({\mathbb R}^n)}
\end{equation}
holds for all $\varphi\in L^2({\mathbb R}^n)$.
\end{lem}
By scaling the proof of (5.7) can be reduced to 
the case $T=1$. 
For $T=1$ the estimate (5.7) 
has been shown in \cite{H2} as a direct consequence of 
integrability (in time) of the local energy \cite{SS}
\begin{equation}
\|W\varphi\|
_{L^2({\mathbb R}\times\{x\in{\mathbb R}^n||x|<1\})}
\leq
C\|\varphi\|_{L^2({\mathbb R}^n)}, 
\end{equation}
scaling, and the energy estimate. 

We are in a position to complete the proof of 
Theorem 5.1. 
Fix any $p$ $(0<1/p\leq 1/2)$ satisfying 
$1/2>(n-1)((1/2)-(1/p))$. 
It follows from 
the Sobolev-type estimate (3.9) and 
(5.7) that 
\begin{eqnarray}
& &
\|W\varphi\|_{L^2((0,T);L^p({\mathbb R}^n))}\\
& &
\leq
C\||x|^{-(n-1)((1/2)-(1/p))}
|D_x|^{(1/2)-(1/p)}W\varphi\|
_{L^2((0,T)\times{\mathbb R}^n)}\nonumber\\
& &
\leq
CT^{(1/2)-(n-1)((1/2)-(1/p))}
\||D_x|^{(1/2)-(1/p)}\varphi\|_{L^2({\mathbb R}^n)}.\nonumber
\end{eqnarray}
Our estimate (5.3) is a consequence of 
the interpolation between (5.9) and 
the energy estimate. 
We have finished the proof of Theorem 5.1.
$\hfill\square$
%%%%%%%%%%%%%%%%%%%%%%%%%%%%%%%%%%%%%%%%
%%%%%%%Section 6%%%%%%%%%%%%%%%%%%%%%%%%
%%%%%%%%%%%%%%%%%%%%%%%%%%%%%%%%%%%%%%%%
\section{End-point estimates for Schr\"odinger equations}
The final section is devoted to the study of 
the Strichartz estimate for the Schr\"odinger equation
\begin{equation}
i\partial_t u-\Delta u=0
\quad
\mbox{in}\,\,{\mathbb R}\times{\mathbb R}^n
\end{equation}
subject to the initial data $u(0,x)=\varphi(x)$. 
The estimate
\begin{equation}
\|S\varphi\|
_{L^2({\mathbb R};L^{2n/(n-2)}({\mathbb R}^n))}
\leq
C\|\varphi\|_{L^2({\mathbb R}^n)}
\quad(n\geq 3),
\end{equation}
which was proved by Keel and Tao \cite{KT}, 
is called an end-point estimate. 
(See (1.2) for the definition of the operator $S$.) 
As Vilela has explained in Section 3 of \cite{V}, 
it is possible to prove (6.2) 
for radially symmetric data 
via the weighted inequality (2.2) with $\alpha=\beta=0$. 
We revisit the problem 
of showing (6.2) for radially symmetric data. 
Using our weighted inequality (2.2) with 
$-\beta=\alpha$, 
we prove
\begin{thm}
Suppose $n\geq 3$ and 
$-(1/2)+(1/n)<\alpha<(1/2)-(1/n)$. 
There exists a constant $C$ 
depending on $n$, $\alpha$, 
and the estimate 
\begin{equation}
\||x|^\alpha |D_x|^\alpha S\varphi\|
_{L^2({\mathbb R};L^{2n/(n-2)}({\mathbb R}^n))}
\leq
C\|\varphi\|_{L^2({\mathbb R}^n)}
\end{equation}
holds for radially symmetric data 
$\varphi\in L^2({\mathbb R}^n)$.
\end{thm}
{\it Proof of Theorem 6.1.}
We need the following lemma. 
\begin{lem}
Suppose $n\geq 2$ and $1/2<\gamma<n/2$. 
There exists a constant $C$ depending on $n$, $\gamma$, 
and the estimate 
\begin{equation}
\||x|^{-\gamma}|D_x|^{1-\gamma} S\varphi\|
_{L^2({\mathbb R}\times{\mathbb R}^n)}
\leq
C\|\varphi\|_{L^2({\mathbb R}^n)}
\end{equation}
holds. 
\end{lem}
Large part of Lemma 6.2 was proved 
by Kato and Yajima \cite{KY}, 
Ben-Artzi and Klainerman \cite{BA-K}, independently. 
Later their results were not only complemented but also generalized 
by Sugimoto \cite{Su} and Vilela \cite{V}. 
For the proof of (6.4) see Section 4 of Sugimoto \cite{Su} 
or Section 1 of Vilela \cite{V}. 

In what follows we denote $2n/(n-2)$ by $p_0$. 
We note that 
$$
\frac12<(n-1)\biggl(\frac12-\frac{1}{p_0}\biggr)-\alpha<\frac{n}{2}
\Longleftrightarrow
1-\frac{1}{n}-\frac{n}{2}<\alpha<\frac12-\frac1n
$$
and that the inequality 
$$
1-\frac{1}{n}-\frac{n}{2}\leq -\frac12+\frac1n
$$
is true for all $n\geq 3$ (actually, for all $n\geq 1$). 
Taking account of the obvious fact
$$
-\alpha<\frac{1}{p_0}
\Longleftrightarrow
-\frac12+\frac1n<\alpha,
$$
we can employing (2.2) with $-\beta=\alpha$ first and (6.4) secondly 
to have for radially symmetric $\varphi$
\begin{eqnarray}
& &
\||x|^\alpha |D_x|^\alpha S\varphi\|
_{L^2({\mathbb R};L^{p_0}({\mathbb R}^n))}\\
& &
\leq
C\||x|^{-(n-1)((1/2)-(1/p_0))+\alpha}
|D_x|^{\alpha+(1/2)-(1/p_0)}\varphi\|
_{L^2({\mathbb R}\times{\mathbb R}^n)}\nonumber\\
& &
\leq
C\|\varphi\|_{L^2({\mathbb R}^n)}.\nonumber
\end{eqnarray}
It is only at the last inequality above that 
the choice of $p_0=2n/(n-2)$ is essential. 
The proof of Theorem 6.1 has been completed.
$\hfill\square$

\vspace{0.3cm}

%%%%%%%%%%%%%%%%%%%%%%%%%%%%%%%%%%%%%%%%%%%%%%%%%%%%%%%
{\bf Acknowledgements.} 
The authors are grateful to Professors Rentaro Agemi and 
Hiroyuki Takamura for their comments. 
They also thank the referee for reading the manuscript carefully and 
making a number of invaluable suggestions. 
The first author was partly supported by the 
Grant-in-Aid for Young Scientists (B) 
(No.\,15740092 and 18740069), 
The Ministry of Education, Culture, 
Sports, Science and Technology, Japan. 
The second author was partly supported by the 
Grant-in-Aid for Young Scientists (B) 
(No.\,18740096), 
The Ministry of Education, Culture, 
Sports, Science and Technology, Japan.
%%%%%%%%%%%%%%%%%%%%%%%%%%%%%%%%%%%%%%%%%%%%%%%%%%%%%%%%

%%%%%%%%%%%%%%%%
% bibliography
%%%%%%%%%%%%%%%

% Set bibliography items using the "thebibliography" environment  and following
% the style used by the AMS journals. 
%
% If the bibliography is generated by a bibtex database, use "amsplain" or
% "amsalpha" as bibliography style

%%%%%%%%%%%%%%%%%%%%%%%%%%%%%%%%%%%%%%%
%%%%%%%%%%%%%%%%%%%%%%%%%%%%%%%%%%%%%%%%

\end{document}